# Global Sensitivity Analysis in Load Modeling via Low-rank Tensor

You Lin, *Student Member, IEEE*, Yishen Wang, *Member, IEEE*, Jianhui Wang, *Senior Member, IEEE*, Siqi Wang, *Member, IEEE,* Di Shi, *Senior Member, IEEE*

*Abstract*—Growing model complexities in load modeling have created high dimensionality in parameter estimations, and thereby substantially increasing associated computational costs. In this paper, a tensor-based method is proposed for identifying composite load modeling (CLM) parameters and for conducting a global sensitivity analysis. Tensor format and Fokker-Planck equations are used to estimate the power output response of CLM in the context of simultaneously varying parameters under their full parameter distribution ranges. The proposed tensor structured is shown as effective for tackling high-dimensional parameter estimation and for improving computational performances in load modeling through global sensitivity analysis.

*Index Terms*--Dimensionality reduction, load modeling, parameter estimation, sensitivity analysis, tensor.

## I. Introduction

LOAD modeling has been a critical component in power system stability analysis for decades. Accurate load modeling leads to more precise system operation limits and can improve system operation and economics [1]. Inaccurate load modeling can result in system operating states that are vulnerable to contingencies, thus emphasizing the need for informative load modeling in power system operations.

Advances in power system measurement equipment, specifically PMUs, allows for the collection of key electrical information of specific buses, such as voltage magnitude $V$, voltage angle $\theta$, real power $P$, and reactive power $Q$. The objective of load modeling, thus, is to construct a mathematical model $f$ to replicate the behaviors of $P$ and $Q$ with input $V$ and $\theta$, as shown in (1).

$$P = f_1(V, \theta|\boldsymbol{\pi}) \text{ and } Q = f_2(V, \theta|\boldsymbol{\pi}) \quad (1)$$

By estimating the load modeling parameters $\boldsymbol{\pi}$, $\boldsymbol{\pi} = [\boldsymbol{\pi}_1, \cdots, \boldsymbol{\pi}_M]^T$, we can calculate the power output response given the input $V$ and $\theta$. Here, $M$ is the total number of parameters, which is usually 4 in ZIP and 11 in composite load modeling (CLM) represented by ZIP and 3rd-order induction motor (IM).

The more complex and detailed the load model structure that is to be mimicked for different types of dynamic responses, the more parameters there are to be estimated. The tradeoff between computation and model structure, therefore, drives the non-trivial nature of the problem for solving and analysis.

This work is funded by SGCC Science and Technology Program under contract no. SGSDYT00FCJS1700676. Y. Lin, Y. Wang, S. Wang and D. Shi are with GEIRI North America, San Jose, CA 95134, USA. Y. Lin and J. Wang are also with the Department of Electrical and Computer Engineering, Southern Methodist University, Dallas, TX 75205. USA.

Additionally, the correlation among parameters increases the difficulty in identifying some parameters. Existing parameter estimation methods in load modeling are highly dependent on the initialization performance.

The existence of multiple locally optimal solutions is a common problem in the current load modeling research field [1]. One common method to obtain a global optimal solution is by calculating the response curve fitting accuracy based on enumerations of all possible parameter sets. However, it is impractical to attain an optimal solution using this method when the number of parameters is large. Moreover, there are several issues still to be addressed in load modeling parameter estimation, including: How sensitive the real and reactive power are to the variation of each parameter? How to efficiently obtain global optimal parameters in their full range?

In this paper, we use tensor structure methods to address the high dimensional parameter estimation problem. We also evaluate the influence of jointly varying parameters to provide global sensitivity analysis. Tensors allow for the mapping of high-dimensional features to low-dimensional representations by explicitly building the tensor in a low-rank form [2]. With such a structure, we can directly control the effective dimensionality of a set of parameters, thereby improving the overall computation performance and solution quality. The contributions of this paper are as follows:
a) Efficient sampling of a large number of parameters in high dimensional parameter space via tensor;
b) Computing varied distributions of all parameters within their entire feasible range;
c) The ability to provide a global optimal reference for all model parameters.

## II. Optimization Problem in Load Modeling

The composite load model in this paper includes the ZIP and IM. The static load component is represented by the ZIP model, whose real power $P_{ZIP,t}$ and reactive power $P_{ZIP,t}$ at time $t$ are calculated from (2) and (3) given the parameters $a_p, b_p, a_q$, and $b_q$.

$$P_{ZIP,t} = P_{ZIP,0}\left(a_p\left(\frac{V_t}{V_0}\right)^2 + b_p\left(\frac{V_t}{V_0}\right) + 1 - a_p - b_p\right) \quad (2)$$

$$Q_{ZIP,t} = Q_{ZIP,0}\left(a_q\left(\frac{V_t}{V_0}\right)^2 + b_q\left(\frac{V_t}{V_0}\right) + 1 - a_q - b_q\right) \quad (3)$$

The dynamic load component is represented by IM. With the measured bus voltage $V_t\angle\theta_t$ and power output $\tilde{S}_t = P_t + jQ_t$, we can derive the d/q transformation of the voltage, $U_{d,t}$ and $U_{q,t}$, from the reference axis of the overall system to d/q axis of motors. The reference axis for the measurements is the global reference axis for the system, which needs to be converted to the d/q axis.

$$\delta_t = tan^{-1}\left(\frac{X_s I_t cos\varphi_t - r_s I_t sin\varphi_t}{V_t - (r_s I_t cos\varphi_t + X_s I_t sin\varphi_t)}\right) \quad (4)$$

where, $\varphi_t = \beta_t - \alpha_t$, and $\dot{I}_t = \left(\frac{\tilde{S}}{V_t \angle \theta_t}\right)^* = I_t \angle \alpha_t$, $r_s$ and $X_s$ are the stator resistance and reactance of IM, respectively. We then can get $U_{d,t} = -V_t \sin(\delta_t)$, and $U_{q,t} = V_t \cos(\delta_t)$.

The third-order IM model can be described by the following differential and algebraic equations (DAEs) [3].

$$\frac{dv'_{d,t}}{dt} = \frac{-r_r}{X_r + X_m}\left(v'_{d,t} + \frac{X_m^2}{X_r + X_m}i_{q,t}\right) + sv'_{q,t} \quad (5)$$

$$\frac{dv'_{q,t}}{dt} = \frac{-r_r}{X_r + X_m}\left(v'_{q,t} - \frac{X_m^2}{X_r + X_m}i_{d,t}\right) - sv'_{d,t} \quad (6)$$

$$\frac{ds_t}{dt} = \frac{1}{2H}\left(T_{m0}(1-s_t)^2 - v'_{d,t}i_{d,t} - v'_{q,t}i_{q,t}\right) \quad (7)$$

$$i_{d,t} = \frac{r_s(U_{d,t} - v'_{d,t}) + X'(U_{q,t} - v'_{q,t})}{r_s^2 + X'^2} \quad (8)$$

$$i_{q,t} = \frac{r_s(U_{q,t} - v'_{q,t}) - X'(U_{d,t} - v'_{d,t})}{r_s^2 + X'^2} \quad (9)$$

$$X' = X_s + \frac{X_m X_r}{X_m + X_r} \quad (10)$$

where $r_r, X_r, X_m, H, s, v'_{d,t}, v'_{q,t}, i_{d,t}, i_{q,t}$ are the rotor resistance and reactance, the excitation reactance, IM inertia and slip, the transient voltages and currents of IM in d/q axis, respectively.

The real and reactive power output, $P_{IM,t}$ and $Q_{IM,t}$, of the IM can be calculated from the following equations:

$$P_{IM,t} = U_{d,t}i_{d,t} + U_{q,t}i_{q,t} \quad (11)$$

$$Q_{IM,t} = U_{d,t}i_{q,t} - U_{q,t}i_{d,t} \quad (12)$$

The estimated real and reactive power output, $\hat{P}$ and $\hat{Q}$, of the composite load model can be calculated from (13) given the static load proportion $\omega$ represented by ZIP.

$$\begin{cases} \hat{P}_t = \omega \cdot P_{ZIP,t} + (1-\omega) \cdot P_{IM,t} \\ \hat{Q}_t = \omega \cdot Q_{ZIP,t} + (1-\omega) \cdot Q_{IM,t} \end{cases} \quad (13)$$

The load modeling formulation in (2)~(13) can be represented with the DAEs as:

$$\begin{cases} d\boldsymbol{x}_t = \boldsymbol{\mu}(\boldsymbol{x}_t, \boldsymbol{\pi})dt \\ \boldsymbol{M}(\boldsymbol{x}_t, \boldsymbol{\pi})\boldsymbol{x}_t = \boldsymbol{b} \end{cases} \quad (14)$$

The optimization problem in load modeling parameter identification is to find the best solution from all the feasible parameter sets $\boldsymbol{\pi}$ and state variables $\boldsymbol{x}$ to achieve the minimized fitting errors. Here, $\boldsymbol{x}_t = \{v'_{d,t}, v'_{q,t}, s_t\}$ and $\boldsymbol{\pi} = \{a_p, b_p, a_q, b_q, r_s, X_s, r_r, X_r, X_m, H, \omega\}$. This fitting error is measured by the root mean square error (RMSE) between the estimated and measured values of real or reactive power.

The load modeling parameter identification is a highly nonlinear and nonconvex optimization problem, which commonly leads to local optimal solutions. These local optimal solutions are also highly dependent on the selection of the initial values. Considering these issues and disregarding the computation limitations, a simple and straightforward method would be to obtain the global optimal solution through enumeration. By calculating the model output accuracy corresponding to all the possible combinations of all dependent parameters in their varying ranges, we can arrive at a global optimal solution.

III. STOCHASTIC SOLUTIONS OF DAE EQUATIONS

Before analyzing the sensitivity of all parameters, we must first find the solution of the DAEs in (15) and then present reasonable distributions of the parameters, which can be achieved simultaneously using the Fokker-Planck operator.

The Fokker-Planck operator is a partial differential equation that describes the time evolution of the joint probability density function of its state variables and parameters [4]. With the DAE equations in (14), the Fokker-Planck equation for the probability distribution $p(\boldsymbol{x}|\boldsymbol{\pi})$ of the state vector $\boldsymbol{x}$ given the parameter values $\boldsymbol{\pi}$ is:

$$\frac{\partial}{\partial t}p(\boldsymbol{x}|\boldsymbol{\pi}) = -\sum_{i=1}^{N}\frac{\partial}{\partial x_i}[\mu_i(\boldsymbol{x},\boldsymbol{\pi})p(\boldsymbol{x}|\boldsymbol{\pi})] \quad (15)$$

where $N$ is the number of state variables. With the parametric Fokker-Planck operator $A(\boldsymbol{x},\boldsymbol{\pi})$, (15) can be expressed as

$$A(\boldsymbol{x},\boldsymbol{\pi})\, p(\boldsymbol{x}|\boldsymbol{\pi}) = -\sum_{i=1}^{N}\frac{\partial}{\partial x_i}[\mu_i(\boldsymbol{x},\boldsymbol{\pi})p(\boldsymbol{x}|\boldsymbol{\pi})] \quad (16)$$

We then use the tensor structures to compute $p(\boldsymbol{x}|\boldsymbol{\pi})$ simultaneously for all parameters $\boldsymbol{\pi}$ within their feasible varying space as shown in the following section.

IV. PARAMETER ESTIMATION VIA TENSOR

In this paper, a tensor structure is used to tackle the high dimension computation burden. Define the $N$-dimensional varying space $\Omega_x = I_1 \times \cdots \times I_N$ of state variables $\boldsymbol{x}_t$ and $M$-dimensional varying space $\Omega_\pi = J_1 \times \cdots \times J_M$ of model parameters $\boldsymbol{\pi}$. Here, the feasible range of state variable $x_{i,t}$ ($i = 1, \cdots, N$) is $I_i = [l_i^x, u_i^x]$, and the feasible range of parameter $\pi_j$ ($j = 1, \cdots, M$) is $J_j = [l_j^\pi, u_j^\pi]$. Here, we discretize the space of $\Omega_x$ into $(n_d + 1)^N$ nodes with step sizes $h_i^x = (u_i^x - l_i^x)/(n_d + 1)$. Similarly, discretizing the space of $\Omega_\pi$ into $(m_d + 1)^M$ points with step sizes $h_i^\pi = (u_i^\pi - l_i^\pi)/(m_d + 1)$. The values of $p(\boldsymbol{x}_t|\boldsymbol{\pi})$ are discretized into $(N+M)$-dimension tensor $\tilde{p} \in R^{i_1 \times \cdots \times i_N \times j_1 \times \cdots \times j_N}$ with entries

$$\tilde{\boldsymbol{p}}_{i_1 \times \cdots \times i_N \times j_1 \times \cdots \times j_N} = p(x_{i_1,t}, \cdots, x_{i_N,t}|\pi_{j_1,t}, \cdots, x_{j_M,t}) \quad (17)$$

If we process (17) in conventional vector and matrix format, the searching space we face is $(n_d + 1)^N (m_d + 1)^M$. The storage cost of tensor $\tilde{p}$ will increase exponentially with $N$ and $M$, which is computationally intractable. Thus, we propose to use the tensor decomposition to approximate the probability distribution $p(\boldsymbol{x}_t|\boldsymbol{\pi})$ [5], shown as:

$$p(\boldsymbol{x}_t|\boldsymbol{\pi}) = \sum_{r=1}^{R} p_1^r(x_{1,t}) \cdots p_N^r(x_{N,t}) p_1^r(\pi_1) \cdots p_M^r(\pi_M) \quad (18)$$

where, $R$ represents the tensor rank. Based on the low-rank representation in (18), the mathematical operations on the probability distribution $p(x_t|\pi)$ in $N+M$ dimensions can be performed using combinations of one-dimensional operations, and the storage cost is reduced to $(N+M)R$.

If choosing a sufficiently large computation domain, i.e., $n_d$ and $m_d$ are large enough, we can approximate the distribution $p(x_t|\pi)$ by eigenvector $\tilde{p}$ of $A(x_t,\pi)$ corresponding to the minimum eigenvalue $\lambda_{min}$ closest to zero. Thus, given the initial value of all state variables $x_0 = \{v'_{d0}, v'_{q0}, s_0\}$ and parameters $\pi_0 = \{a_{p0}, b_{p0}, a_{q0}, b_{q0}, r_{s0}, X_{s0}, r_{r0}, X_{r0}, X_{m0}, H_0, \omega_0\}$, $\tilde{p}$ can be solved from (19) using the alternating minimum energy method (AMEN) method [6] as the linear solver in tensor format.

$$A(x_t,\pi)\tilde{p} = \lambda_{min}\tilde{p} \quad (19)$$

Since the solution $\tilde{p}$ is presented in tensor format shown in (18), we extract the univariate probability density function of each state variable and parameter, represented by $p_i(x_{i,t})$ and b $p_j(\pi_j)$, respectively.

By examining the probability densities of each parameter, the best parameter estimation $\pi^*$ at its largest probability is obtained. Simultaneously, several local optimal estimations at its lower probability positions are also provided.

Given the optimal parameter estimation $\pi^*$, the real and reactive power output of the CLM can be obtained using (2)~(13). Meanwhile, the model fitting error can be calculated.

## V. FRAMEWORK OF THE PROPOSED MODEL

The pseudo code of the tensor-based parameter estimation algorithm is presented in Table I.

TABLE I. THE TENSOR-BASED PARAMETER ESTIMATION ALGORITHM

| Pseudo code of the tensor-based parameter estimation algorithm |
|---|
| **Input:** Measured bus voltage magnitude $V_t$ and voltage angle $\theta_t$ |
| **Initialize** the state variables $x_0 = \{v'_{d0}, v'_{q0}, s_0\}$ and model parameters $\pi_0 = \{a_{p0}, b_{p0}, a_{q0}, b_{q0}, r_{s0}, X_{s0}, r_{r0}, X_{r0}, X_{m0}, H_0, \omega_0\}$. |
| Define the varying range of $x_t$ and $\pi$. |
| Formulate the DAEs $dx_t = \mu(x_t,\pi)dt$. |
| Calculate the Fokker–Planck operator $A(x_t,\pi)$ in tensor format. |
| while $A(x_t,\pi)p(x_t|\pi) \neq 0$ |
|    Adjust the AMEN algorithm. |
|    $i=i+1$. |
| **Output:** Individual distributions of $x_t$ and $\pi$, represented by $p_i(x_{i,t})$ and $p_j(\pi_j)$; global optimal solution for model parameters $\pi^*$; and response of power output $\hat{P}$ or $\hat{Q}$ corresponding to $\pi^*$. |

## VI. NUMERICAL RESULTS

The real parameters used in the simulation are shown in Table II. The distribution of each parameter is obtained from the formulated model. The estimated joint distributions of motor inertia $H$ and motor rotor reactance $X_r$ are shown in Fig. 1. The estimated distributions of $H$ and $X_r$ concentrate in values around 0.95 and 0.25, respectively. Similarly, the estimated joint distributions of static load proportion $\omega$ and motor stator reactance $X_s$ are shown in Fig. 2. The estimated optimal values of $\omega$ and $X_s$ are 0.4974 and 0.1053, respectively. From the distributions in Fig. 1 and 2, local optimal values of these two parameters are also found. The model power responses (P/Q) resulting from the voltage is further obtained by combining all distributions of parameters and state variables.

## VII. CONCLUSION

A global sensitivity analysis approach is proposed to perform load modeling parameter estimation via tensor. Simultaneously, varying distributions of all load model parameters in their feasible range can be efficiently estimated based on the Fokker–Planck equations in tensor format. The low-rank tensor representation of the formulated load modeling can tackle high dimensional problems with efficient computation. Global optimal references of parameters are obtained from the estimated distribution of all load model parameters. Furthermore, the sensitivities of all parameters are intuitively presented in the estimated parameter distributions.

TABLE II. REAL MODEL PARAMETERS

| Parameter | $\omega_1$ | $a_p$ | $b_p$ | $a_q$ | $b_q$ | |
|---|---|---|---|---|---|---|
| Value | 0.5 | 0.001 | 0.5642 | 0.001 | 0.6626 | |
| Parameter | $r_s$ | $X_s$ | $r_r$ | $X_r$ | $X_m$ | $H$ |
| Value | 0.049 | 0.096 | 0.044 | 0.244 | 2.96 | 0.93 |

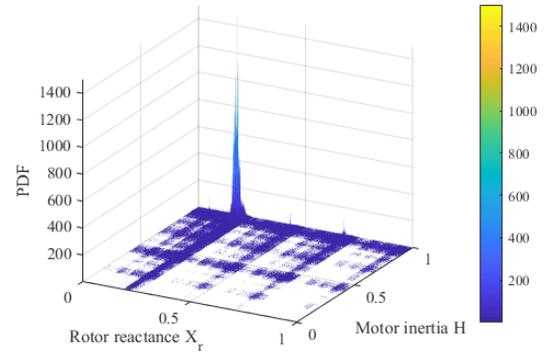

Fig.1 Simulated PDF of model parameters ($H$ and $X_r$)

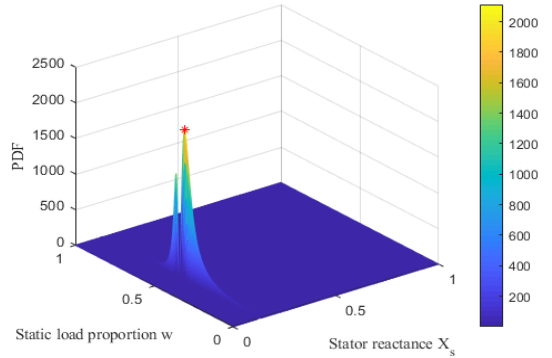

Fig.2 Simulated PDF of model parameters ($\omega$ and $X_s$).